\documentclass[referee,pdflatex,sn-mathphys-num]{sn-jnl}

\usepackage{graphicx}%
\usepackage{multirow}%
\usepackage{amsmath,amssymb,amsfonts}%
\usepackage{amsthm}%
\usepackage{mathrsfs}%
\usepackage[title]{appendix}%
\usepackage{xcolor}%
\usepackage{textcomp}%
\usepackage{manyfoot}%
\usepackage{booktabs}%
\usepackage{algorithm}%
\usepackage{algorithmicx}%
\usepackage{algpseudocode}%
\usepackage{listings}%

\theoremstyle{thmstyleone}%
\newtheorem{theorem}{Theorem}
\newtheorem{conjecture}{Conjecture}%
\newtheorem{open question}{Open question}%

\theoremstyle{thmstyletwo}%

\theoremstyle{thmstylethree}%

\begin{document}

\title[Article Title]{Growth Forms of Tilings}


\author[1]{\fnm{Peter} \sur{Hilgers}}\email{peter.hilgers@posteo.de}

\author[2]{\fnm{Anton } \sur{Shutov}}\email{a1981@mail.ru}

\affil[1]{ \orgaddress{\street{Hasenbergweg 9}, \city{Ettlingen}, \postcode{76275}, \country{Germany}}}

\affil[2]{\orgname{Vladimir State University}, \orgaddress{\street{Gorky Str. 87}, \city{Vladimir}, \postcode{600000}, \country{Russian Federation}}}



\maketitle

\section{Introduction}\label{sec1}
\renewcommand{\figurename}{Figure}

Tilings have been an aspect of human culture since prehistoric times, and the standard reference for their mathematical treatment is~\cite{Grünbaum1987}. In this paper we focus on one feature of tilings, the so-called growth forms.
In our opinion, growth forms are a fascinating but neglected topic in tiling theory. The growth forms of periodic tilings are completely understood; however, for cut-and-projection tilings, general solutions are missing~\cite{Akiyama2019}, and our knowledge is limited to the case of special windows~\cite{Demski2022}. In crystallography coordination sequences~\cite{Brunner1993} or topological densities~\cite{OKeeffe1991}, geometric or topological characteristics, are well established. Growth forms help to study the former~\cite{Shutov2019CoordPer,Shutov2020CoordZuni} or to compute the latter~\cite{Eon2004TopDens,Shutov2020TopPer,Shutov2025TopDens}. Growth forms have also been used to study NaCl clusters, sulphur nanoclusters, carbon single crystals, or the density of states of a Schr\"odinger operator ~\cite{Rau2002,Rau2009,Khorkov2018,Azamov2022}. Structures closely related to growth forms are used in~\cite{Kotani2006LargeDev} to study random walks on a crystal lattice.

The purpose of this paper is to collect the scattered results that have been obtained in recent years \cite{Akiyama2016,Akiyama2019,Demski2022,Lutfalla2025}, to report on the results of computer experiments, and to illustrate all the results with often surprising and aesthetic examples.
For our proofs, we either refer to the literature, or provide them in an appendix.

A tiling is a covering of a $d$-dimensional space by figures with nonoverlapping interiors.
The figures forming a tiling are called tiles.
We say that two tiles are neighbours if they have an intersection of positive $(d - 1)$-dimensional volume.
For example, for $d = 2$, two tiles have common parts of their boundaries of positive length.

The concept of a coordination shell originates in crystallography and is defined inductively:
A patch $P_0$ is a finite set of tiles in the tiling.
The first coordination shell $P_1$ of $P_0$ consists of all tiles that are adjacent to a tile of size $P_0$.
The n-th coordination shell $P_n$ contains all tiles adjacent to $P_{n-1}$, which are not in $P_{n-2}$.
The union of the first $N$ coordination shells is known as the $N$-corona.
The use of the terms 'neighbour' and 'corona' is not homogeneous in the literature. Discussions about different definitions of coronas can be found in~\cite{Dolbinin2000}.
For example, studies of Heesch's tiling problem use the requirement to have an intersection of positive $(d - 2)$-dimensional volume as a neighbourhood definition.
The corona definition changes accordingly.
We use the described definition for two reasons: it is the traditional definition used in all previous work on the subject, and it has a crystallographic motivation. If two-dimensional tiles have the same edge, this can be interpreted as a strong interaction between molecules in a crystal.

\begin{figure}[ht]
\centering
\includegraphics[width=0.9\textwidth]{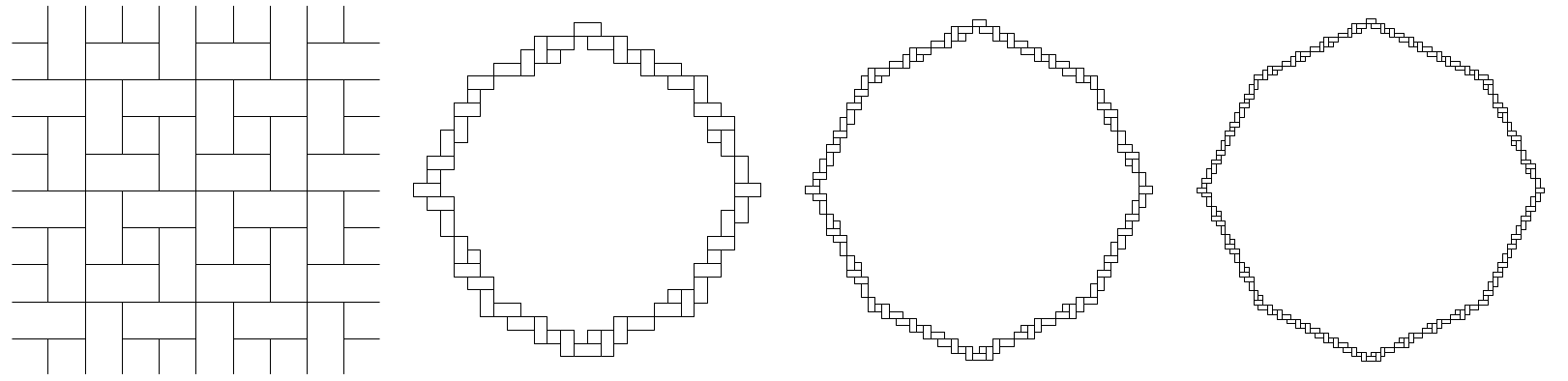}
\caption{Tiling, scaled coordination shells $P_8, P_{16}, P_{24}$. }\label{fig:scaled shells}
\end{figure}

If we scale the coordination shells $P_n$ by the factor $1/n$, we can see that in many cases the obtained sets tend to some limited form, e.g., Figure~\ref{fig:scaled shells}.
This form is called the growth form of tiling.
Informally, we can write $form = lim\;P_n/n$.
More formally, the closed $(d-1)$-dimensional surface form is called a growth form of the tiling if a point $x$ and a constant $C$ exist so that the $n$-th coordination shell lies in the $C$-neighbourhood of the set $x + n*form$.
Further specification of this definition is possible using the concept of the Hausdorff distance.
If the growth form exists, it does not depend on the initial patch $P_0$~\cite{Zhuravlev2002}.

\section{Periodic tilings}\label{periodic}

Periodic tilings can be described using lattices:
A d-dimensional lattice $L$ is a set of vectors of the form $n_1v_1 +\ldots+ n_dv_d$ where the vectors $v_1,\ldots, v_d$ lie in $d$-dimensional Euclidean space and are linearly independent and the numbers $n_1,\ldots, n_d$ are integers.
Two tiles are called equivalent modulo lattice if one of them can be obtained from the other by translation through some vector of $L$.

A tiling is periodic if some $L$ and a finite number of tiles, called fundamental tiles, exist so that any other tile of the tiling is equivalent to one of the fundamental tiles modulo lattice.
Note that, generally, there are many ways to choose fundamental tiles.
The union of fundamental tiles is known as the fundamental domain.
\begin{theorem}
\label{theorem periodic}
Every periodic tiling has a growth form which is a convex centrosymmetric polytope.
\end{theorem}

Different proofs of Theorem~\ref{theorem periodic} can be found in~\cite{Rau2002,Zhuravlev2002,Fritz2013,Akiyama2019}, with similar results in~\cite{Inoue2024Erhard}. Algorithm~\ref{algorithm periodic}~\cite{Zhuravlev2002,Fritz2013} facilitates the computation of the growth form of periodic tilings. The characteristics of all possible growth forms are described in~\cite{Shutov2014Inverse}.

\begin{algorithm}
\caption{Calculate the growth form of a periodic tiling.}\label{algorithm periodic}
\begin{algorithmic}[0]
\Require Periodic tiling $T$
\Ensure Growth form
\State \verb+Find set+ $F$\verb+ of fundamental tiles+
\State $z \Leftarrow \verb+ Number of tiles in + F$
\State $L \Leftarrow \verb+ Lattice corresponding to +T$
\ForAll{$f \in F$}
	\State $shells \Leftarrow\verb+ First +z\verb+ coordination shells+$
	\ForAll{$s \in shells$}
		\State $\verb+Find tiles +$T'$ \verb+ equivalent to +T\verb+ mod +L$
		\State $\verb+Denote by +$v$ \verb+ the vector fullfilling +T'=T+v$
		\State $\verb+Assume that +T'\verb+ lies in the +$k$\verb+-th shell, then + $w$=v/k $
		\State $\verb+Collect all +w$
	\EndFor	
\EndFor
\State$\verb+The growth form is the convext hull of all +w.$
\end{algorithmic}
\end{algorithm}
The reader may attempt to find the growth form of the regular tiling $4^4$ ~\cite{Grünbaum1987} of the plane with squares of area 1 using Algorithm~\ref{algorithm periodic}. It is a square with vertices $v_1 = (1, 0), v_2 = (0, 1), v_3 = (-1, 0), v_4 = (0, -1)$.

\begin{figure}[H]
\centering
\includegraphics[width=0.9\textwidth]{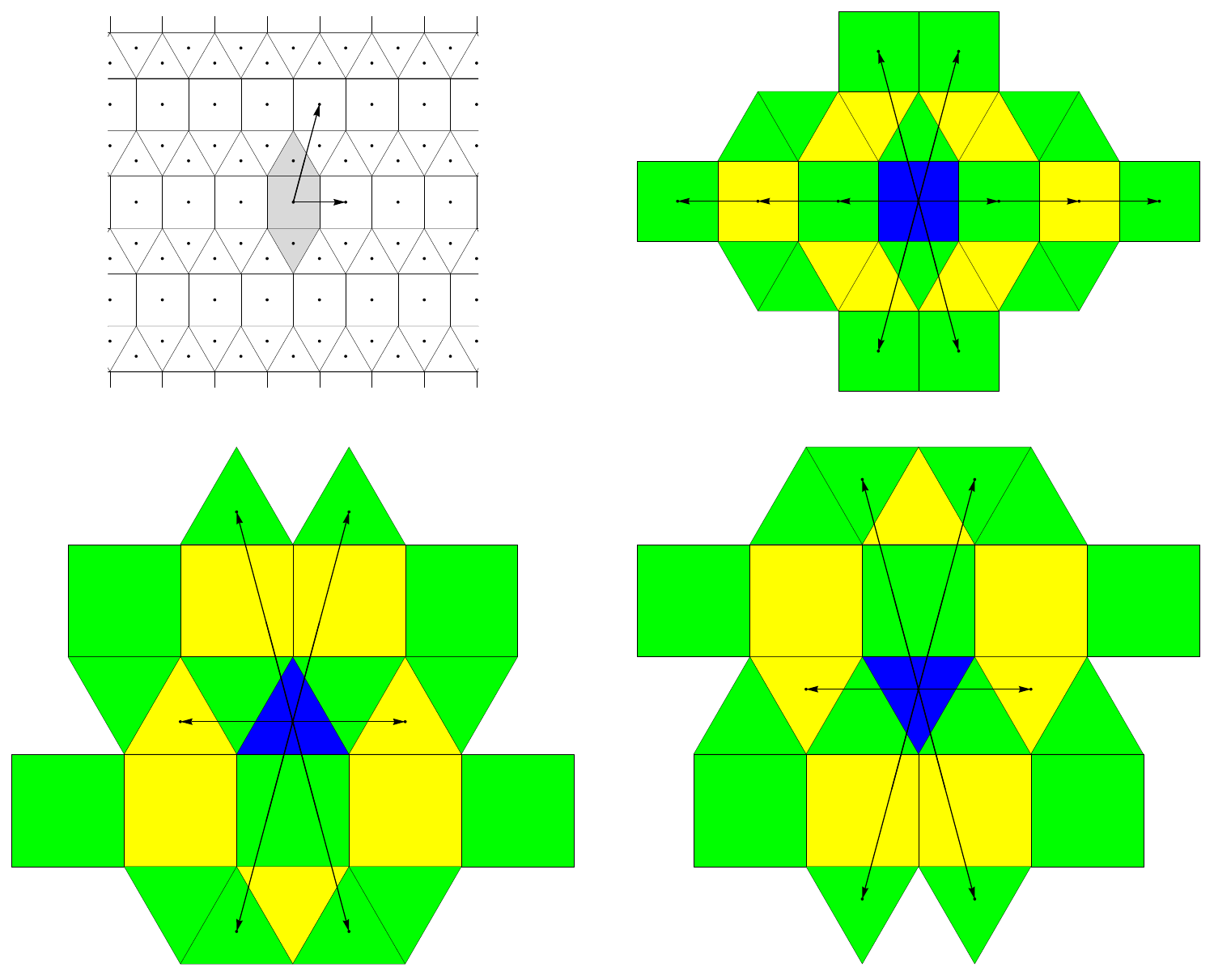}
\caption{Top left: Tiling, fundamental tiles (grey), lattice, and translation vectors. Others: Fundamental tiles (blue) with $z = 3$ coordination shells and vectors $v$.}
\label{fig periodic 1}
\end{figure}

\begin{figure}[H]
\centering
\includegraphics[width=0.9\textwidth]{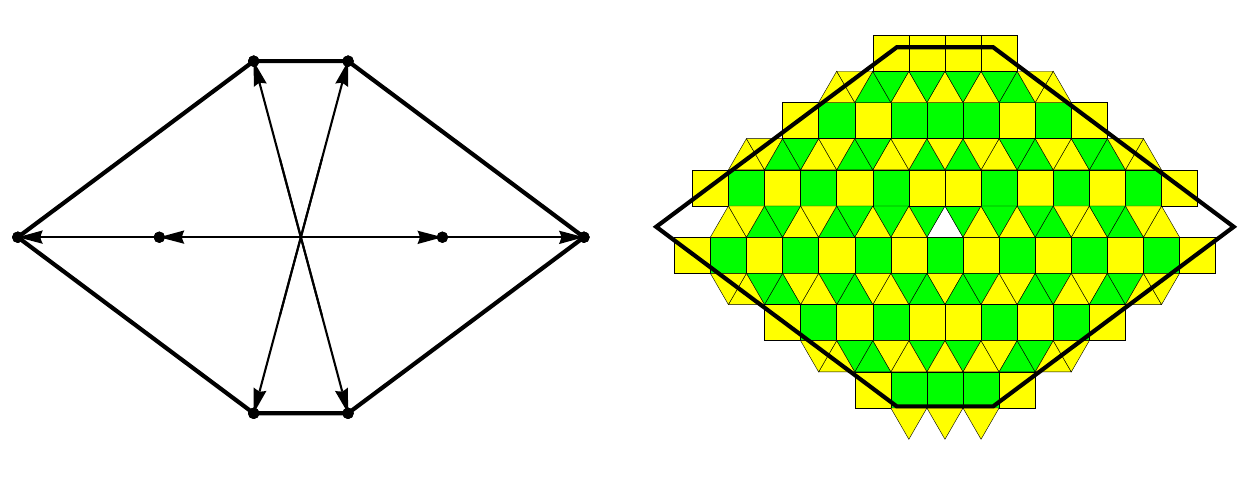}
\caption{Left: Vectors $w$ and growth form. Right: Growth form, patch $P_0$ (white), 12-corona scaled, odd coordination shells green, even ones yellow.}
\label{fig periodic 2}
\end{figure}

As an example, we present one of the eleven Archimedean tilings, the $(3^3.4^2)$ tiling~\cite{Grünbaum1987} with three fundamental tiles, Figures~\Ref{fig periodic 1}, \Ref{fig periodic 2}.


Note that near some vertices of the proven growth form, the speed of convergence seems to be slow; this can also be observed in Figures~\ref{fig grid examples} and~\ref{fig growth form hat}.

\section{Grid tilings in two and three dimensions}\label{grid}

In 1981, de Bruijn~\citep{Bruijn1981} devised an algebraic description of rhombic Penrose tilings on the basis of the dualisation of a pentagrid.
Variations and generalisations of the method are used to study quasicrystals~\cite{Kramer1984,Madison2020,Boyle2022}. These methods are essentially equivalent to the so-called cut-and-projections methods~\cite{Gaehler1986}.

Let ${\{g_1,\ldots, g_N\}}$ be a family of $N$ unit vectors in $d$-dimensional Euclidean space.
We also choose $N$ real parameters $\gamma_i$, which serve as phase shifts.
The $N$-grid $L_N$ is a union of $N$ arrays of equidistant parallel hyperplanes in $\mathbb{R}^d$:
\begin{equation}
L_N=\{x\in \mathbb{R}^d:(x,g_i)-\gamma_i\in \mathbb{Z}
\end{equation}

Here, $(\bullet,\bullet)$ is a scalar product, and $1 \le i \le N$.
If there is no point where more than $d$ grid hyperplanes intersect, the grid is called \textbf{regular}.
Note that the grid will be regular for almost all values of $\gamma_i$.
The $N$-grid $L_N$ defines some tiling of the $d$-dimensional space, but this tiling is ''bad'' because the number of tile types is infinite.
The key idea of the grid method is to consider a tiling dual to $L_N$ in some sense.
Define $N$ functions $K_i$ and function $K$ as follows:

\begin{equation}
K_i(x)=\min\lbrace n\in \mathbb{Z}\ge(x,g_i)-\gamma_i\rbrace
\end{equation}

\begin{equation}
K(x)=\sum_{i=1}^NK_i(x)g_i
\end{equation}

\begin{figure}[H]
\centering
\includegraphics[width=0.9\textwidth]{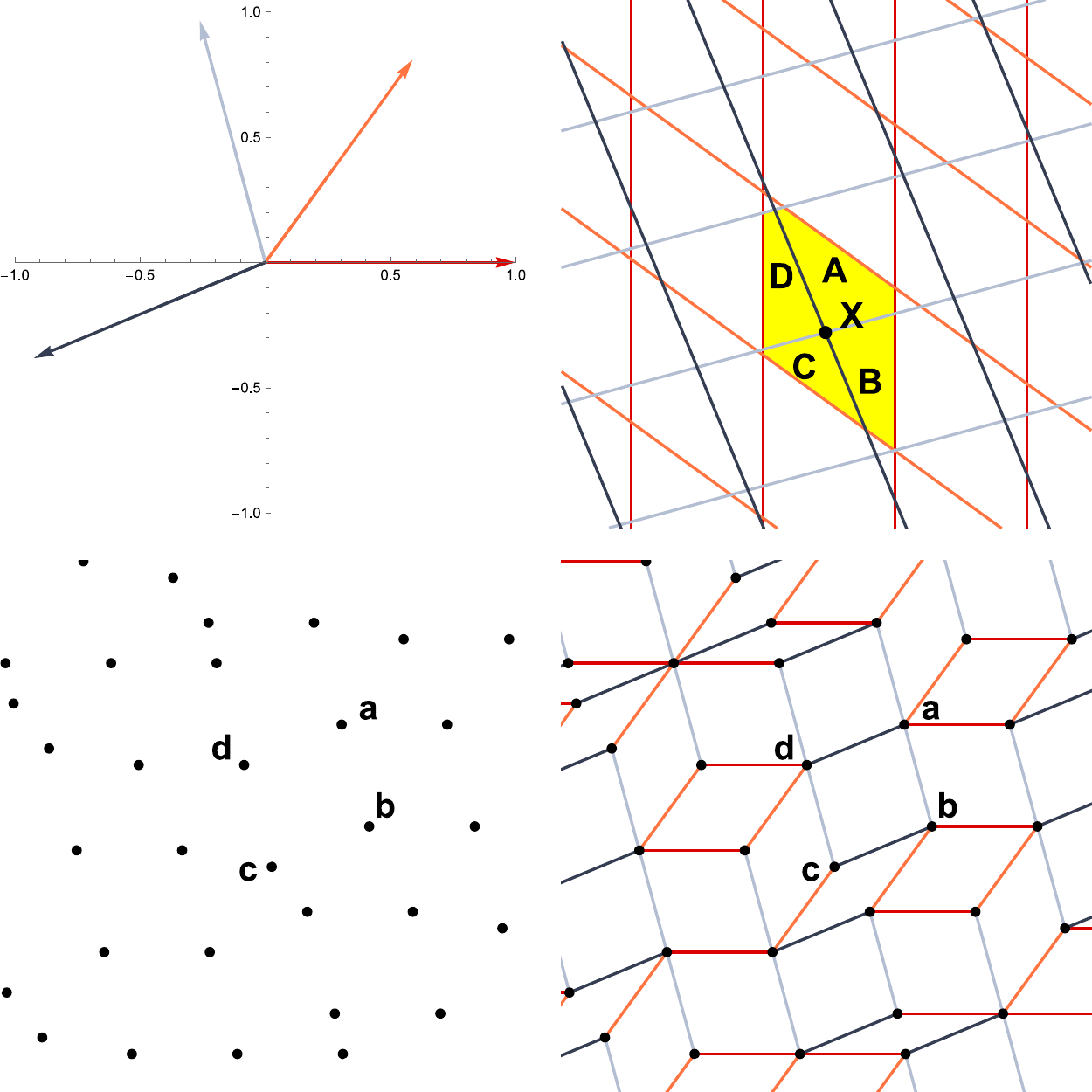}
\caption{Construction of a grid tiling: top left: grid vectors $g_i$; top right: 4-grid $L_N$ with intersection point $X$; adjacent tiles A, B, C, D; bottom left: vertices $\Lambda$ of the dual tiling with A$\leftrightarrow$a, B$\leftrightarrow$b, C$\leftrightarrow$c, D$\leftrightarrow$d; bottom right: edges of dual tiling.}
\label{fig grid constructio}
\end{figure}

Informally, $K_i(x)$ is the index of the hyperplane perpendicular to $g_i$ through $x$.
It can be proven that $K(x)$ is constant on tiles of $L_N$. Therefore, $K$ maps the set of tiles of $L_N$ to a discrete set $\Lambda$ in $\mathbb{R}^d$. The set $\Lambda$ is a set of vertices of some $d$-dimensional tiling, which is called \textbf{grid tiling}.
To define this tiling, one must describe the edges connecting points from $\Lambda$.
The rule is as follows: two points of $\Lambda$ are connected by an edge if and only if the corresponding tiles of $L_N$ have a common edge.
Furthermore, the set of all tiles of $L_N$ sharing some fixed vertex is mapped (under $K$) to the set of all vertices of some tile of the grid tiling.
In the three-dimensional case, the grid tiles are parallelepipeds.
The construction of a two-dimensional tiling is exemplarily illustrated in Figure~\ref{fig grid constructio} for $N = 4$ and $d = 2$.

Now, consider a regular $N$-grid tiling produced by vectors ${\{g_1,\ldots,g_N\}}$.
Let ${\{e_1,\ldots,e_N\}}$ be the standard orthonormal base of the $N$-dimensional space.
Let $O_N$ be a boundary of the convex hull of the vectors $\pm e_i$.
The polytope $O_N$ is known as the $N$-dimensional cross-polytope or as the $N$-dimensional \textbf{orthoplex}.
Define the projections
\begin{align}
&\pi_1:\mathbb{R}^N\rightarrow \mathbb{R}^d \mbox{ defined as }\pi_1(t)=\pi_1((t_1,\dotsc,t_N))=\sum_{i=1}^N t_ig_i\mbox{ and}\\
&\pi_2:\mathbb{R}^N\rightarrow \mathbb{R}^{N-d}\mbox{,}
\end{align}
an orthogonal projection to the $(N-d)$-dimensional plane $\pi_1(t)=0$. Let $P$ be the $d$-dimensional plane $\pi_2(t)=0$.

\begin{theorem}
\label{theorem regular grid}
For any regular grid tiling, the growth form exists and is $\pi_1(O_N \cap P)$~\cite{Demski2022}.
\end{theorem}

\begin{figure}[H]
\centering
\includegraphics[width=0.5\textwidth]{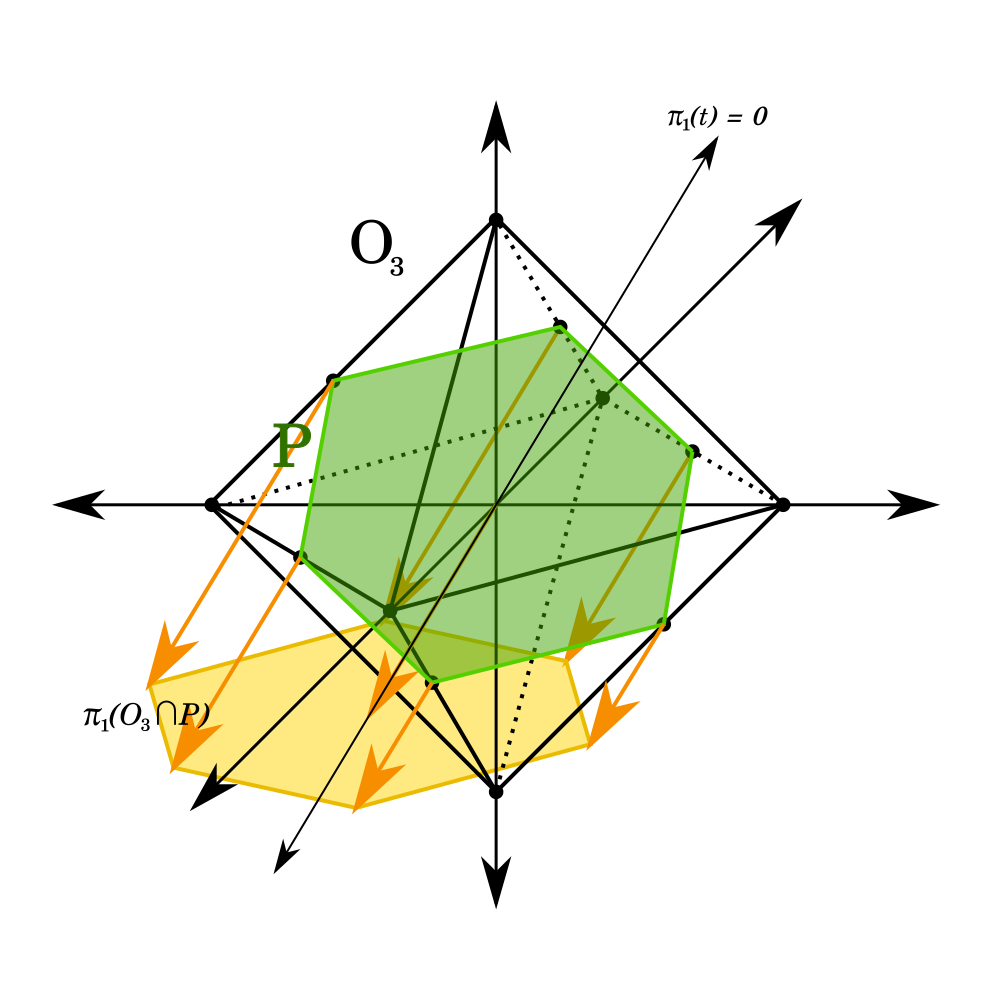}
\caption{: Example of finding the growth form for a case where $N = 3$, $d = 2$}
\label{fig grid find growth form}
\end{figure}

Figure~\ref{fig grid find growth form} illustrates the theorem.

Appendix A explains why it is possible to give explicit formulas for the coordinates of the growth form vertices:
We formulate here the results for the two- and three-dimensional cases:

In the two-dimensional case, let $\alpha_{ij}$ be acute angles between lines orthogonal to grid vectors $g_i$ and $g_j$.
$\alpha_{ij}$ does not depend on the choice of concrete lines, and can be computed as $\alpha_{ij}=\arccos\lvert(g_i,g_j)\rvert$.
Let $\delta_{ij}$ be the distance between two points of intersection of a line orthogonal to $ g_i $ with two consecutive lines orthogonal to $ g_j $.
$\delta_{ij}$ also does not depend on the choice of concrete lines, and can be computed as $\delta_{ij}=1/\sin\alpha_{ij}$.

Define
\begin{align}
D_{ij} &= \begin{vmatrix}g_{i1}&g_{i2}\\g_{j1}&g_{j2}\end{vmatrix}\,,&
\epsilon_{ij}&=\bigg\lbrace\begin{matrix}1,&D_{ij}> 0\\-1,&D_{ij}< 0\end{matrix}\:,&
\Delta_{i}&=\sum_{j\neq i}\frac{1}{\delta_{ij}}\,,&
\upsilon_{i}=\sum_{j\neq i}\frac{\epsilon_{ij}g_j}{\delta_{ij}}\:.
\end{align}

\begin{theorem}
\label{theorem regular 2d grid}
The growth form of a regular two-dimensional grid tiling is a polygon with vertices
\begin{equation}
\left\{\pm\frac{\upsilon_i}{\Delta_i}\right\}, i=1,\ldots,N.
\end{equation}
\end{theorem}

In the three-dimensional case, we use the following notation:
$P_k$ --plane orthogonal to $g_k$, $l_{ij}$ -- line of intersection of planes $P_i$ and $P_j$, $\alpha_{ijk}$ -- angle between $l_{ij}$ and $P_k$, $e_1 = (1, 0, 0)$, $e_2 = (0, 1, 0)$, $e_3 =(0, 0, 1)$, $h_{ij}=g_i\times g_j$.
Note that $h_{i,j}$ is not a unit vector.
Then, $\alpha_{ijk}=\arcsin\frac{(h_{ij},g_{k})}{\lvert h_{ij}\rvert}$. Let $ \delta_{ijk} $ be the distance between two points of intersection of line $l_{ij}$ with two consecutive planes orthogonal to $ g_{k} $.
Then, $ \delta_{ijk}=\frac{1}{\lvert \sin \alpha \rvert}=\frac{\lvert h_{ij} \rvert}{\lvert (h_{ij},g_{k})\rvert} $.

Define
\begin{align}
D_{ijk} &= \begin{vmatrix}g_{i1}&g_{i2}&g_{i3}\\g_{j1}&g_{j2}&g_{j3}\\g_{k1}&g_{k2}&g_{k3}\end{vmatrix}\,,&
\epsilon_{ijk}&=\bigg\lbrace\begin{matrix}1,&D_{ijk}> 0\\-1,&D_{ijk}< 0\end{matrix}\:,&
\Delta_{ij}&=\sum_{k\neq i,j}\frac{1}{\delta_{ijk}}\,,&
\upsilon_{ij}=\sum_{k\neq i,j}\frac{\epsilon_{ijk}g_k}{\delta_{ijk}}\:.
\end{align}

\begin{figure}[H]
\centering
\includegraphics[width=0.9\textwidth]{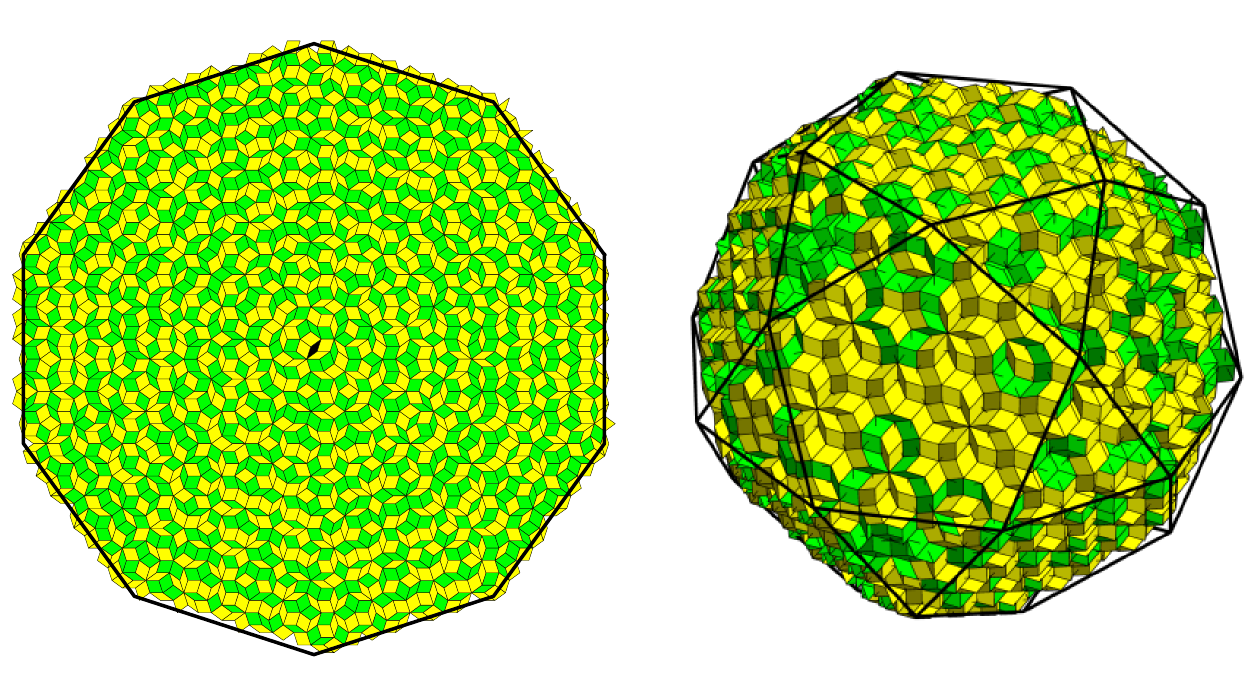}
\caption{Growth forms and scaled n-coronas, left: Penrose tiling, n = 29; right: Ammann 3D tiling, n = 19.}
\label{fig grid examples}
\end{figure}

\begin{theorem}
\label{theorem regular 3d grid}
The growth form of a regular three-dimensional grid tiling is a polyhedron with vertices
\begin{equation}
\left\{\pm\frac{\upsilon_{ij}}{\Delta_{ij}}\right\}, 1\le i<j\le N.
\end{equation}
\end{theorem}
The proof of Theorem~\ref{theorem regular 3d grid} is similar to that of Theorem~\ref{theorem regular 2d grid}. Table~\ref{tableTwoGrowthForms} and
Figure~\ref{fig grid examples} exemplifies the theorems.

\begin{table}[h]
\centering
\begin{tabular}{{@{}lll@{}}}
\toprule
Tiling &Penrose&Ammann 3D~\cite{Senechal2004,Socolar1986} \\
\midrule
Growth form~\cite{Demski2022} & Regular decagon&Icosidodecahedron\\
Circumradius growth form~\cite{Demski2022} &$5/2\tan{2\pi/10}=5/2\sqrt{5-2\sqrt{5}}$&$\sqrt{5-2\sqrt{5}}$\\
Colour coding in Figure \ref{fig grid examples} &Coordination shell&Type of tile\\
\bottomrule
\end{tabular}
\caption{Growth forms of two tilings }
\label{tableTwoGrowthForms}
\end{table}

\section{Substitution tilings}\label{substitution}

A substitution tiling is a tiling that is built by a simple method:
Use a finite set of building blocks, the prototiles $T_1, T_2,\ldots ,T_m$.
Expand all prototiles by a linear map $Q$, the inflation factor.
Apply a set of rules for dissecting each scaled tile $Q\:T_i$ into copies of the original prototiles.
Following this process \textit{ad infinitum} fills the whole space.
The growth form of tilings with an inflation structure is a large, unsolved problem.
There are fundamental questions:
\begin{open question}
\label{open question subst univ}
Is there a universal method to decide whether a substitution tiling has a growth form?
\end{open question}
\begin{open question}
\label{open question subst algo}
If it exists, is there a general method to determine the growth form?
\end{open question}

Even the published understanding of special cases is limited to the case of pinwheel tiling:
It has a circular growth form~\cite{Radin1996}.
Therefore, we restrict our discussion to two remarkable examples. These are the chair and the L-tetramino tiling.~\cite{Frettloeh,
Grünbaum1987}.

\subsection{Chair tiling}\label{chair}
\begin{figure}[H]
\centering
\includegraphics[width=0.9\textwidth]{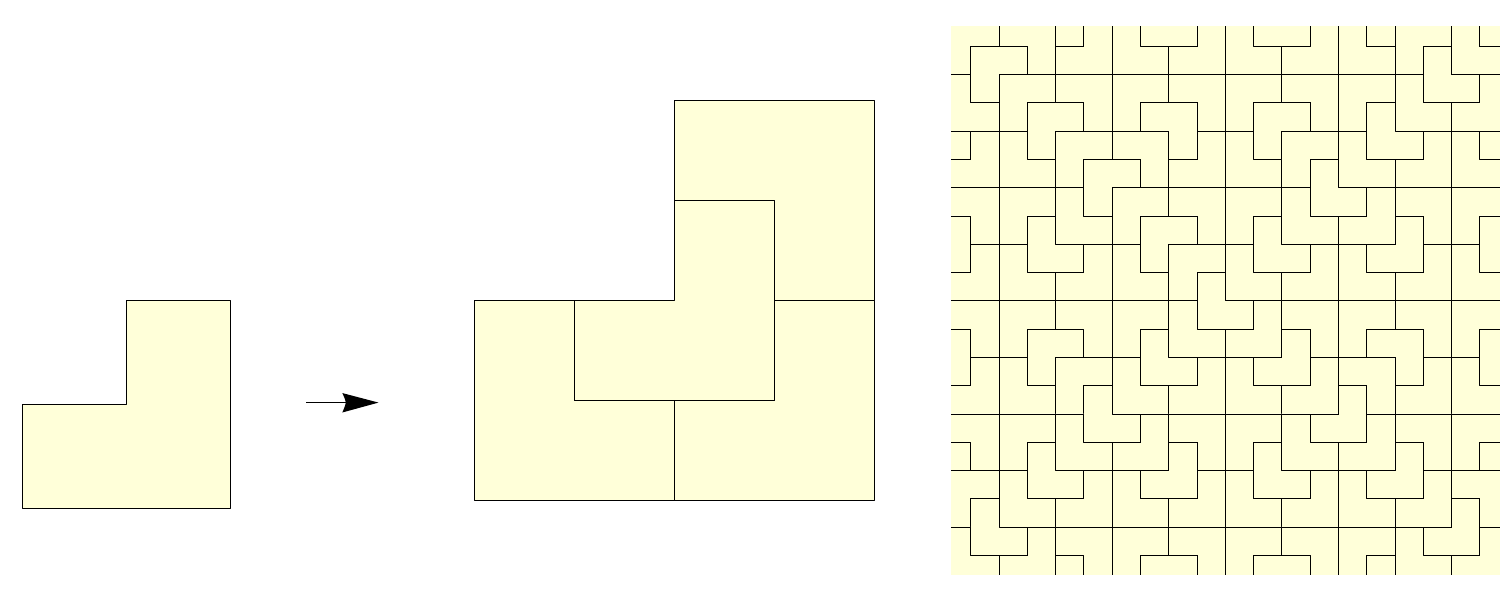}
\caption{Chair tiling. Left: prototiles, inflation and dissection; right: tiling.}
\label{fig chair 1}
\end{figure}

The chair tiling is one of the 3-ominoe tilings~\cite{Grünbaum1987}, Figure~\ref{fig chair 1}.

\begin{figure}[H]
\centering
\includegraphics[width=0.9\textwidth]{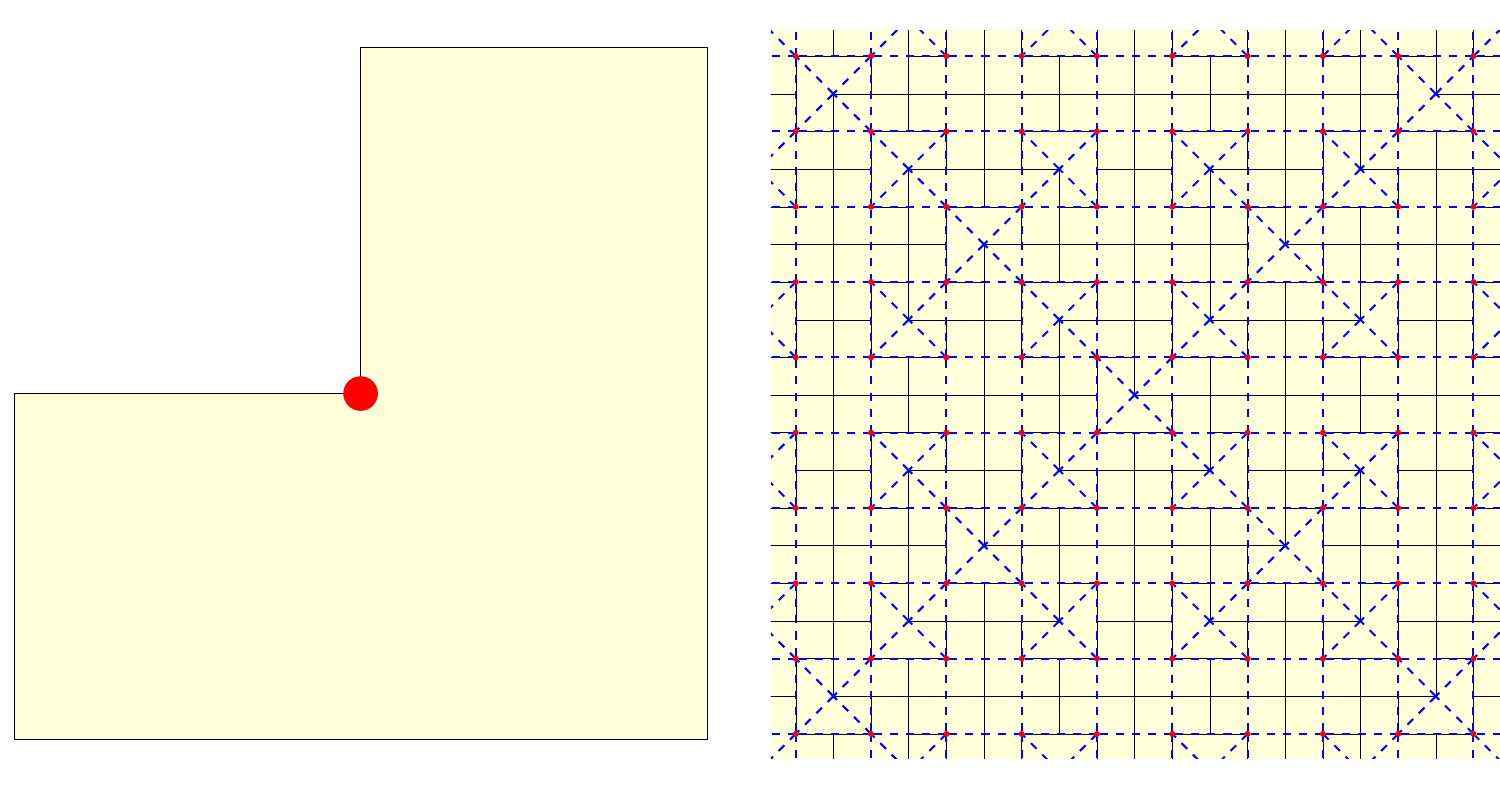}
\caption{Dual tiling to chair tiling. left: prototile; right: tiling and dual tiling.}
\label{fig chair 2}
\end{figure}

In each chair choose a point, as shown in Figure~\ref{fig chair 2}, and construct a dual graph in the following way.
Two points are connected by an edge if and only if the corresponding chairs are neighbours.

\begin{figure}[H]
\centering
\includegraphics[width=0.5\textwidth]{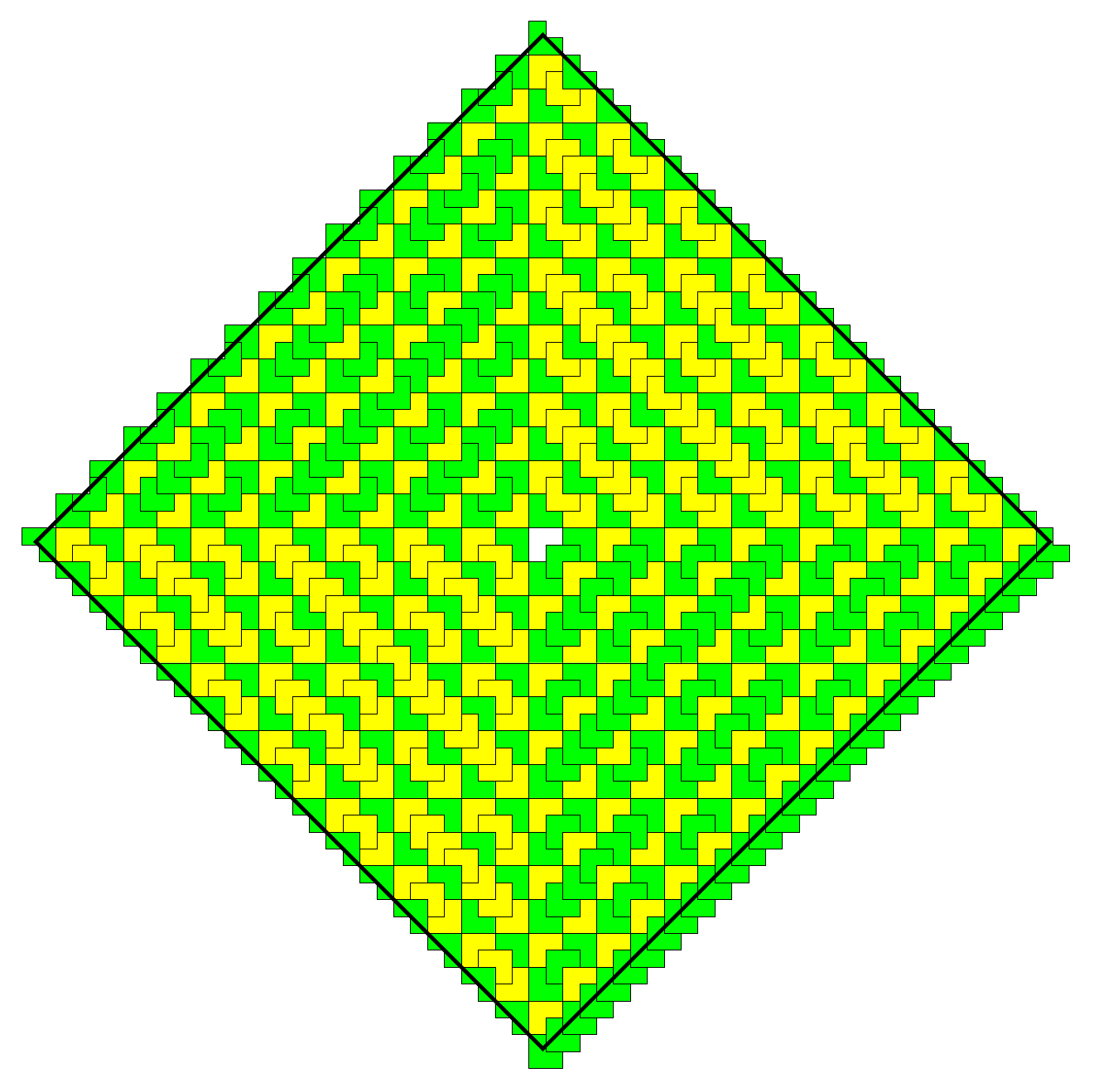}
\caption{Scaled coordination shells $P_1$ to $P_{15}$ of the chair tiling, growth form.}
\label{fig chair 3}
\end{figure}

Therefore, in dual tiling, we have a combination of two structures, squares, and squares with additional diagonals.
It is easy to see that deleting the diagonals does not change the distance between other vertices.
Hence, we have:

\begin{theorem}
\label{theoreme chair}
The growth form of the chair tiling is identical to the growth form of the periodic square graph (Figure~\ref{fig chair 3}).
\end{theorem}

\subsection{L-tetromino tiling}\label{L tetromino}
The L-tetromino tiling has one prototile:
\begin{figure}[H]
\centering
\includegraphics[width=0.9\textwidth]{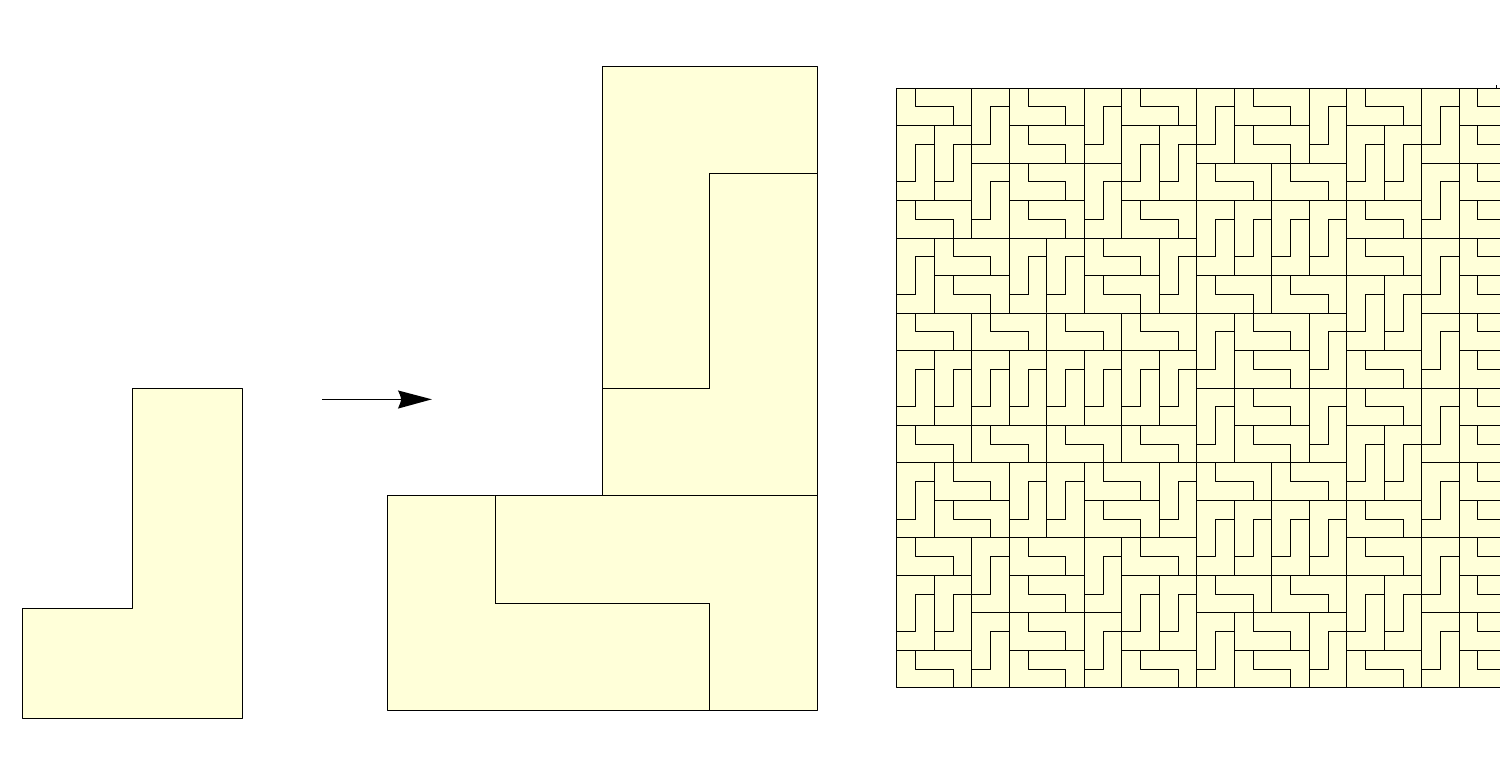}
\caption{L tetromino tiling, left: prototiles, inflation and dissection, right: tiling. }
\label{fig L-tetro 1}
\end{figure}
Our experimental results lead us to suppose the following:
\begin{conjecture}
\label{conjecture L}
The growth form of L. tetromino tiling is not a polygon. Figure~\ref{fig L-tetro 2}.
\end{conjecture}
\begin{figure}[H]
\centering
\includegraphics[width=0.5\textwidth]{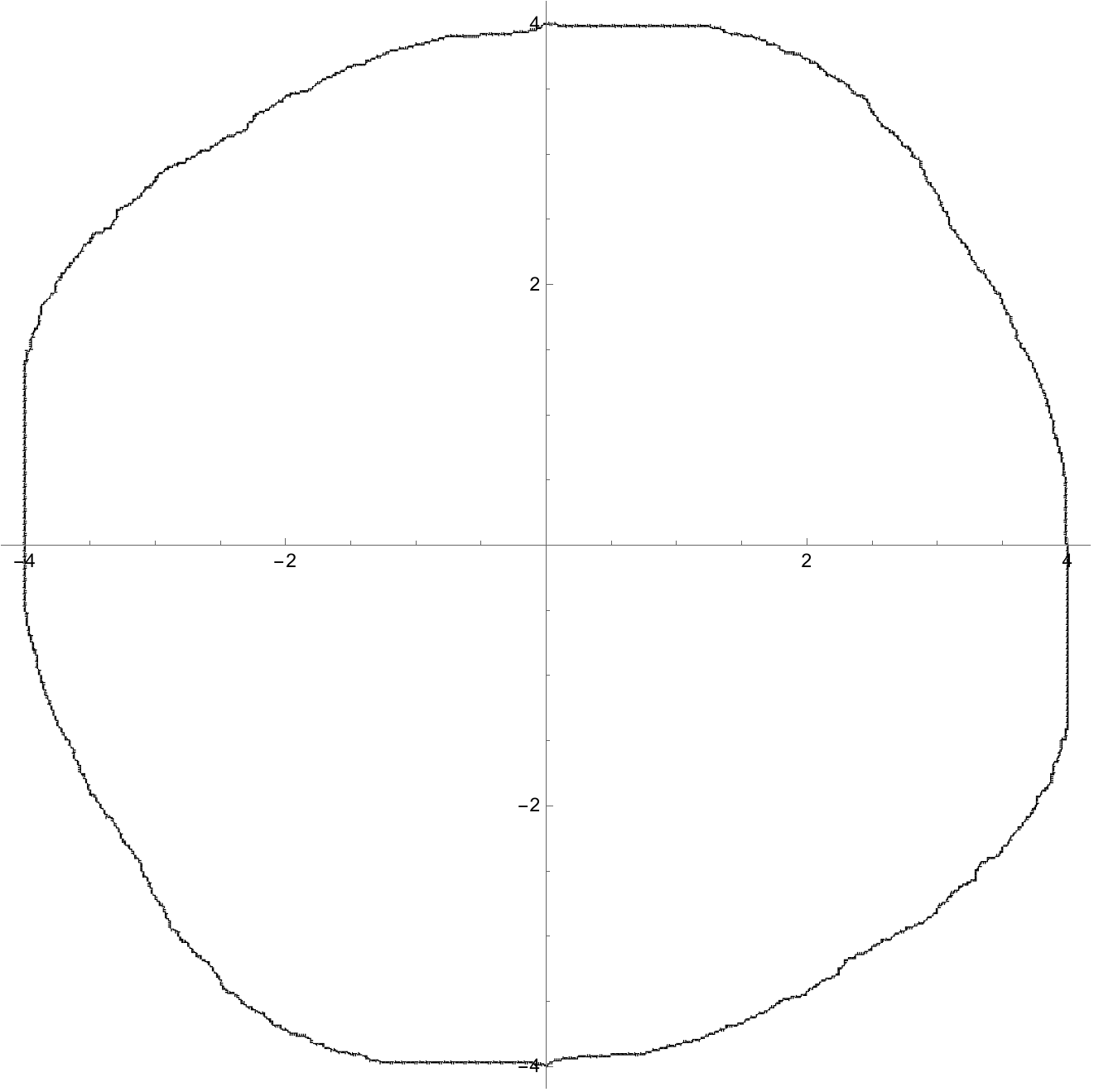}
\caption{L-tetromino tiling, scaled coordination shell $P_{255}$, 1651 tiles.}
\label{fig L-tetro 2}
\end{figure}
Further experimental results can be found in~\cite{Imai2018pdf}.

\section{Hat and related tilings}\label{hat}

\begin{figure}[H]
\centering
\includegraphics[width=0.5\textwidth]{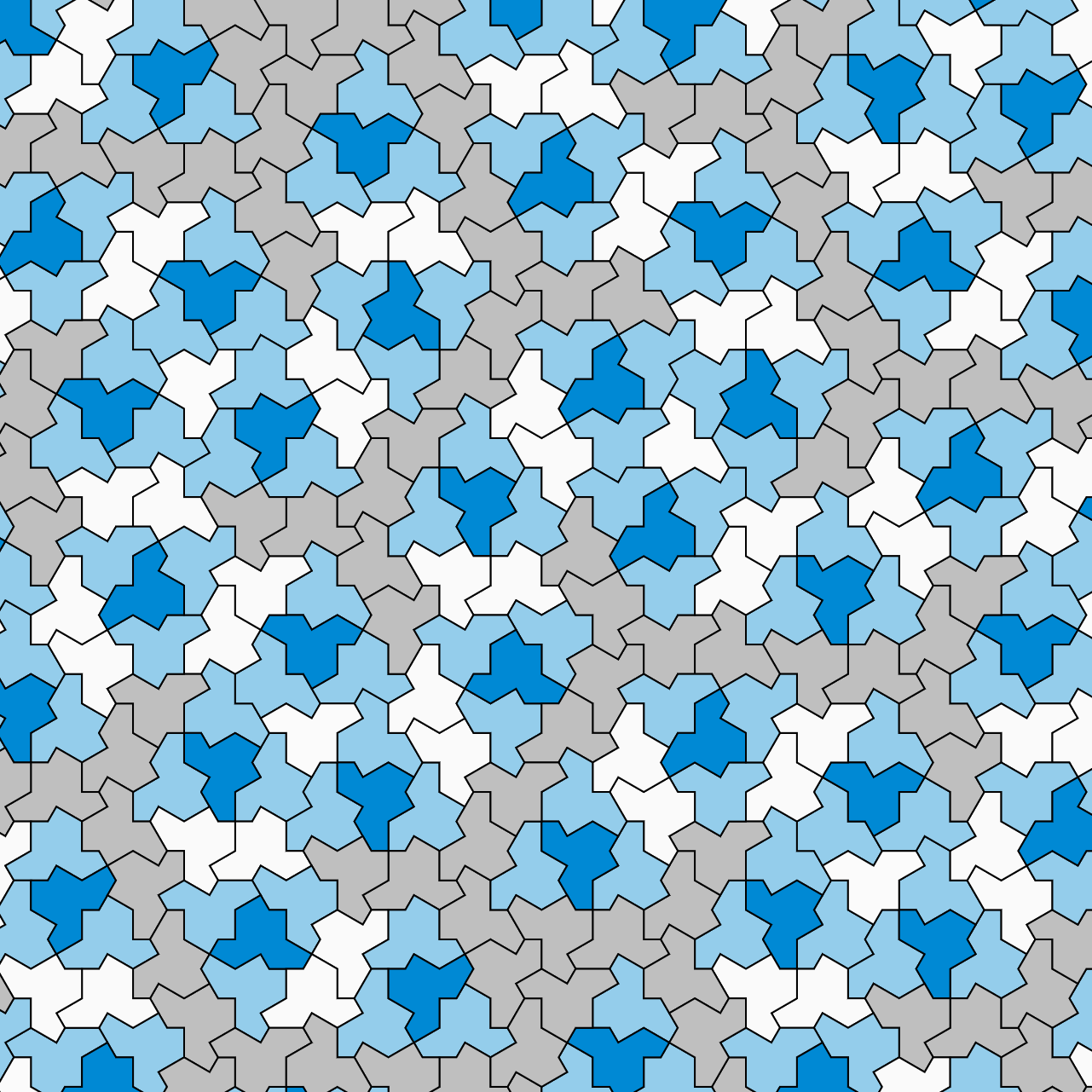}
\caption{Hat tiling~\citep{Smith2024}.}
\label{fig Hat Kaplan coloured}
\end{figure}

Hat tiling is a new fascinating substitution tiling that was discovered as a solution to the famous ''Ein stein'' problem.
There are various alternative constructions of this tiling~\cite{Smith2024,Socolar2023,Reitebuch2023}.
Smith et al.~\cite{Smith2024} reported that hat tiling $HAT$ has a natural colouring by 4 colours, usually denoted as dark blue, light blue, grey, and white (Figure~\ref{fig Hat Kaplan coloured}).

\begin{theorem}
\label{theorem Hat}
Hat tiling has a regular hexagon as growth form.
\end{theorem}

\begin{figure}[H]
\centering
\includegraphics[width=0.9\textwidth]{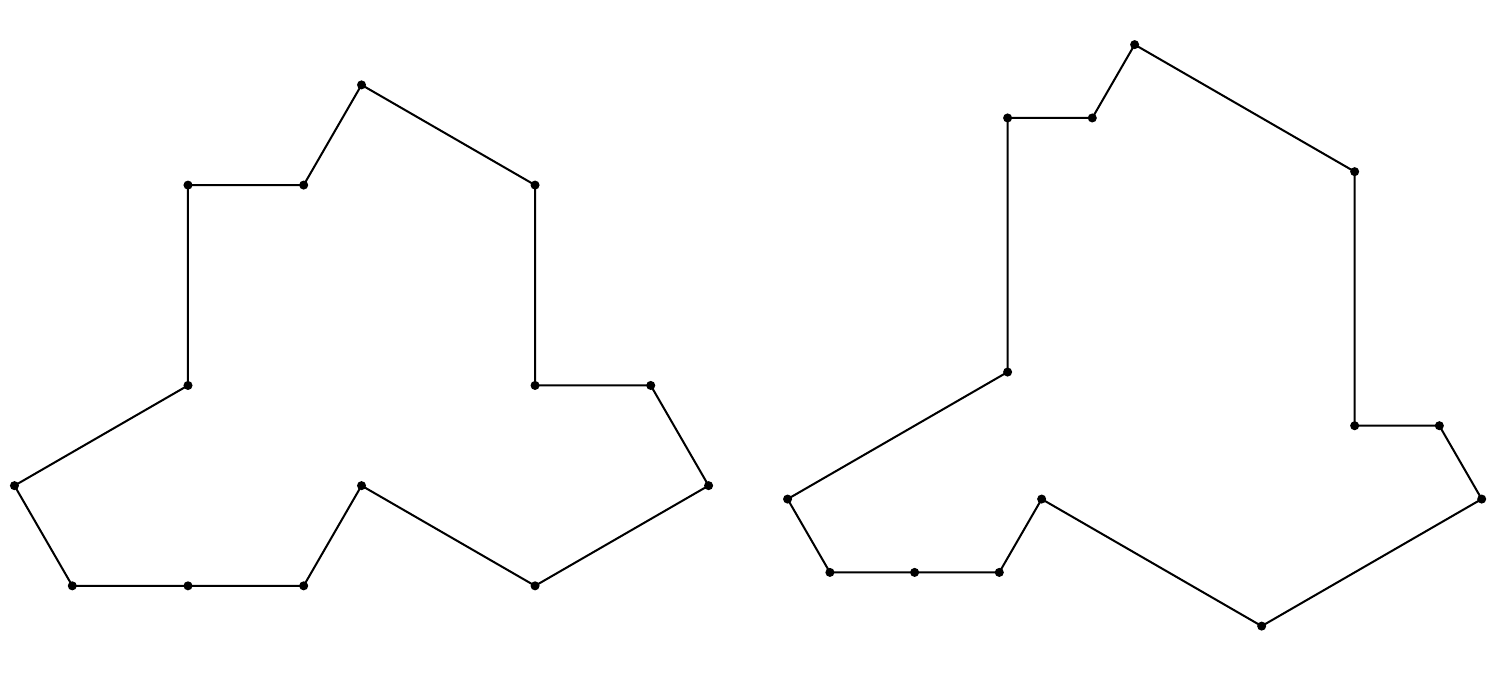}
\caption{Tiles: Hat and $Tile(1, 3)$ }
\label{fig Hat tile forms}
\end{figure}

\begin{figure}[H]
\centering
\includegraphics[width=0.9\textwidth]{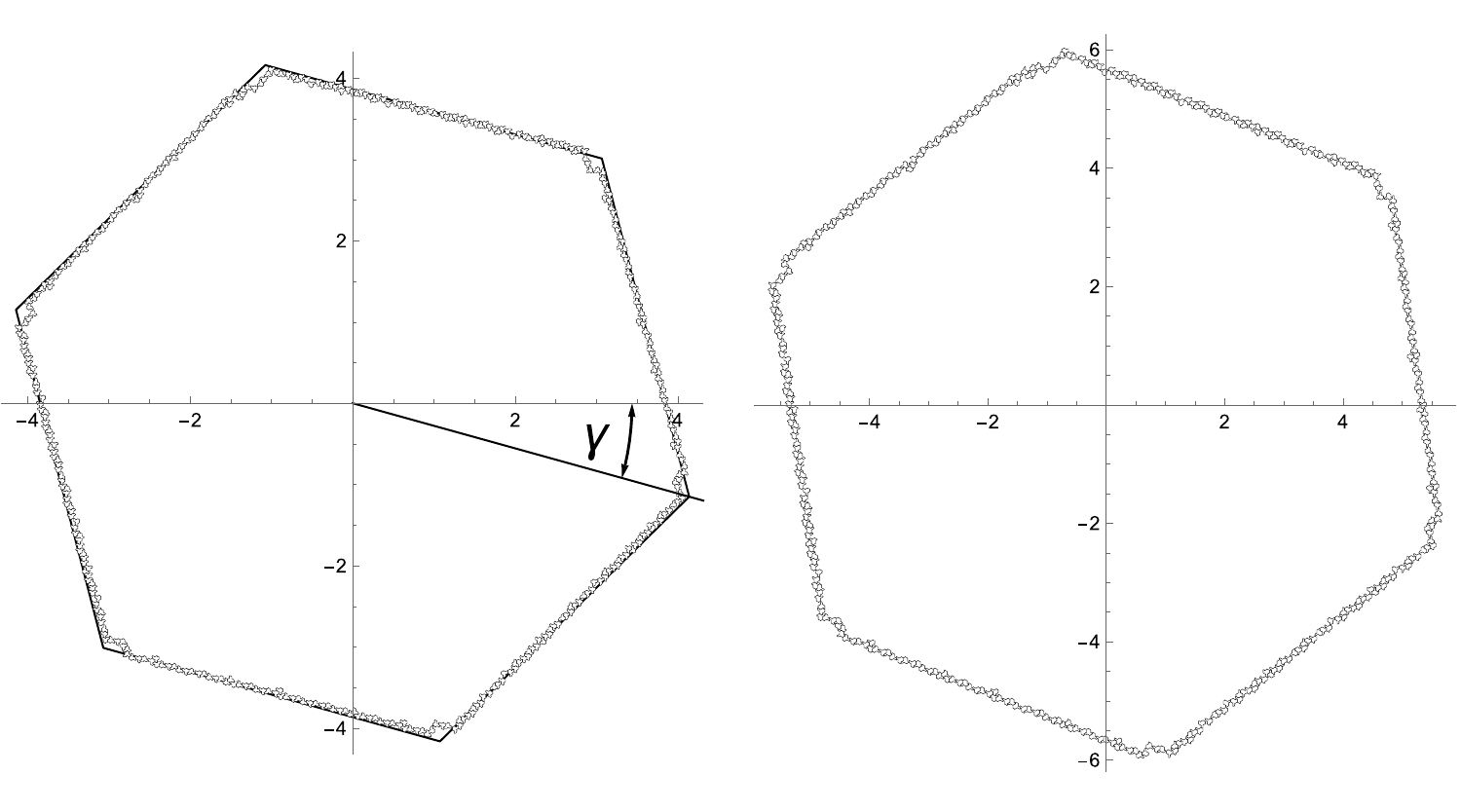}
\caption{Left: Growth from Hat tiling, scaled shell $P_{55}$. Right: scaled shell $P_{55}$ from $Tile(1, 3)$; note the different scales.}
\label{fig growth form hat}
\end{figure}

The colouring is related to the so-called metatile structure, which we do not need, but the colouring is used in the proof of the theorem, Appendix B.
The discoverers of Hat tiling also introduce deformed tiles called $Tile(a, b)$, as shown in Figure~\ref{fig Hat tile forms}.
As tiles $Tile(a, b)$ and $Tile(1, b/a)$ are similar, we can consider only tiles $Tile(1, b)$ with positive $b \neq 1$.
All resulting tilings are aperiodic.
They contain $Tile(1, b )$ in 6 rotational orientations and their mirror images. The hat tile is $Tile(1, \sqrt{3})$ and has edge lengths of 1 and $\sqrt{3}$.
The area of $Tile(1, b)$ is $\sqrt{3} (2+\sqrt{3}b+b^2 )$.
It is also proven that all corresponding tilings are combinatorially equivalent.
Because of this, we can deduce that the growth forms of all such tilings are regular hexagons.

Currently, we do not know the edge lengths and tilt angles of the growth hexagons.
On the basis of computer experiments, and the fact that hat tiling is closely related to the hexagonal lattice, we propose the following:
\begin{conjecture}
\label{conjecture Tile(1,b)}
The growth forms resulting from $Tile(1, b)$ are as follows:
\begin{align*}
area_{growth form}= 2\sqrt{3}\:area_{Tile(1, b)}
\end{align*}
\end{conjecture}

The edge lengths follow immediately, as shown in Table~\ref{tbl:parameter}.

\begin{table}[h]
\centering
\begin{tabular}{{@{}lll@{}}}
\toprule
Tile &Hat&$Tile(1,b)$ \\
\midrule
b&$\sqrt{3}$&3\\
Area tile & $8\sqrt{3}$&$9+11\sqrt{3}$\\
Area growth form & 48&$16+18\sqrt{3}$\\
Edge length growth form &$4\sqrt{2}/\sqrt[4]{3}$&$2\sqrt{3+11/\sqrt{3}}$\\
Tilt angle $\gamma$ &$\arctan(\sqrt{3}/(3+2\tau))$&?\\
\bottomrule
\end{tabular}
\caption{Parameter of tiles and growth forms }
\label{tbl:parameter}
\end{table}

The triangulation lines in Figure~\ref{fig Hat proof 1} have slopes of $\alpha$ and $\alpha+2\pi/6$ with $\tan(\alpha) = \sqrt{3}/(3+2\tau)$, where $\tau$ is the golden ratio~\citep{Reitebuch2023}. We conjecture:

\begin{conjecture}
\label{conjecture Hat tilt}
The growth form of Hat tiling is tilted by $\gamma = -\alpha = -0.270919$.
\end{conjecture}

\section{A tiling with no growth form}\label{no growth form}
\begin{figure}[H]
\centering
\includegraphics[width=0.9\textwidth]{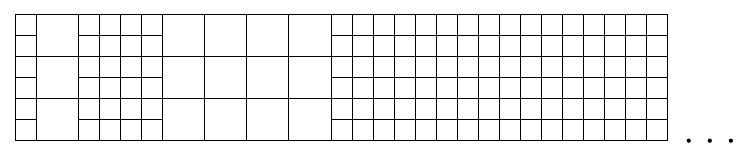}
\caption{A tiling with no growth form. }
\label{fig no growth form}
\end{figure}

Not all tilings have a growth form. A simple example is a tiling by 1 $\times$ 1 and 2 $\times$ 2 square tiles arranged in vertical strips of increasing width, Figure~\ref{fig no growth form}. Let s be a strip of 1 $\times$ 1 tiles, and l be a strip of 2 $\times$ 2 tiles. The sequence of strips is then defined by $4^is\;4^il, i=0\ldots\infty$. For proofs and details, see~\citep{Lutfalla2025}.

\section{Conclusion}\label{sec13}

We have presented the concept of the growth form of tilings and illustrated it with many examples.
For three classes of tilings, namely, the periodic, the regular grid tilings in two and three dimensions, general methods exist to compute the growth forms.
However, for other types of tilings, our knowledge is limited and incomplete.
In the wide range of substitution tilings, we only know the growth form of the pinwheel tiling, which is a circle~\cite{Radin1996}, of the chair and the hat tiling.
The last two cases have not yet been published.
Computer experiments have led to conjectures about, e.g., L-tetromino tiling.

For us, the most important open questions concern substitution tilings: Is there a general method to resolve whether a substitution tiling has a growth form?
Does an algorithm exist to compute the growth form as the case may be?
Which values of the substitution matrix induce the existence of a growth form?
Another important question is the following:
What will change if we use the Heesch definition of  the corona?
In the periodic case, we will have the same algorithm for finding the growth form and the same proof as before.
However, in the nonperiodic case, all the known proofs do not work.
Computer experiments show that if we have a polygonal growth form with our definition of the corona, then the growth form based on the Heesch definition will be another polygon, but there is not a single case of proof.

Special aspects also pose open questions:
Can the growth form be a fractal?
We do not know, but on the basis of~\cite{Shutov2014Inverse}, we conjecture that the answer is no, and that a growth form must be a boundary of some star body.
We also observe that nonregular grid tilings may have nonconvex growth forms.
It is unclear under which conditions convex and nonconvex forms occur.

\begin{appendices}
\section{Sketch of proof of theorem~\ref{theorem regular 2d grid}}\label{Appendix grid}

We define a chain as a sequence of tiles so that any tile in the sequence is neighbouring the previous tile.
The length of the chain is the number of tiles in the chain minus one. It suggests itself to use the minimal length of the chain from tile x to tile y as a distance function of a tiling. Geodesic chains are chains of minimal length between two tiles.
In the two-dimensional case, geodesic chains in the direction of growth from vertices are dual to grid lines.
It can easily be seen that this observation is valid in arbitrary dimensions.

Choose some grid lines orthogonal to $g_i$.
Then, for any segment of this line having length $l$, the number $N_{ij}(l)$ of points of intersection of this line with grid lines orthogonal to $g_j$ satisfies
\begin{equation}
\lim_{l\rightarrow\infty }\frac{N_{ij}(l)}{\delta_{ij}}=1.
\end{equation}
Moreover, some $C$ exist such that $ |N_{ij}-l/\delta_{ij}|\le C $ is independent of the choice of line and segment. Furthermore, if $ N_{i}(l)$ is the number of all points of the multigrid on the segment, we have $ N_{i}(l)=\sum_{j\neq i} N_{ij}(l)$.
Hence,
\begin{equation}
\lim_{l\rightarrow\infty}\frac{N_i(l)}{l/\Delta_i}=1.
\end{equation}
Furthermore, let x be some point of a tile of size $L_N$ neighbouring the beginning of the segment.
Consider a chain of tiles of length $L_N$ having length $ N_i(l) $, with the property that any tile of the chain is neighbouring the chosen line.
Let $y(l)$ be a point from the last tile of the chain.
Because geodesic chains in the direction of growth from vertices are dual to grid lines, vertices corresponding to the grid line orthogonal to $g_i$ are $\pm\lim_{l\rightarrow\infty }\frac{K(y(l))-K(x)}{N_i(l)} $.
We must verify that this limit does not depend on the choice of x or the concrete grid line.
Resulting from the definitions, we have $K_{j}(y(l))-K_{j}(x)=\epsilon_{ij}N_{ij}(l)$.
Hence, $K(y(l))-K(x)=\sum_{j\neq i}\epsilon_{ij}N_{ij}(l)g_{j}  $.
Therefore, using (A2), we obtain
\begin{equation}
\lim_{l\rightarrow\infty}\frac{K(y(l))-K(x)}{l}=\frac{\sum_{j\neq i}\epsilon_{ij}l/\delta_{ij}g_{j}}{l}=\upsilon_{i}
\end{equation}
Finally, using (A1) and combining (A2) and (A3), we obtain the required result.
Note that a result equivalent to Theorem~\ref{theorem regular 2d grid} is also obtained in~\cite{Lutfalla2025}.

\section{Sketch of the proof of theorem \ref{theorem Hat}}\label{secA3}%
\begin{figure}[H]
\centering
\includegraphics[width=0.9\textwidth]{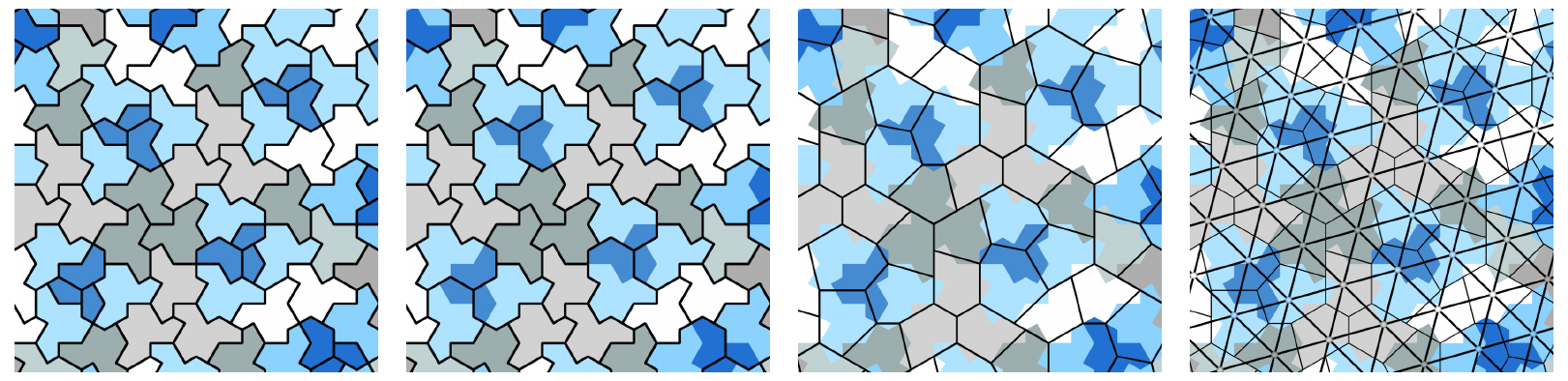}
\caption{Left to right: Split flipped hats into three parts – Join parts to adjacent light blue hats – Hexagon mesh connectivity – Dual triangle mesh~\cite{Reitebuch2023}.
}
\label{fig Hat proof 1}
\end{figure}

Our proof uses the idea of Reitebuch~\cite{Reitebuch2023} to split each dark blue tile into 3 parts.
Each of them will be added to the adjacent light blue tile.
Therefore, we obtain a new plane tiling (coloured by 3 colours).
Denote this as $HAT_1$. This new tiling is a deformation of the tiling of the plane with regular hexagons.
To achieve this goal, we can straighten the sides of the hat tiling to obtain plane tiling by convex polygons. After that, we take a point in each tile and connect points by edges if and only if the corresponding tiles are neighbours. This construction is known as dual tiling. As a result, we obtain a tiling of the plane by regular triangles, as shown in Figure~\ref{fig Hat proof 1}.

This tiling is also dual to the hexagonal tiling of the plane. From this, we can deduce that $HAT_1$ and the hexagonal tiling of the plane have the same growth form. A sketch of the proof is as follows:
\begin{enumerate}
\item Consider the graph growth form of dual triangle tiling. This graph is periodic, and its growth form is a regular hexagon that coincides with the growth form of the hexagonal tiling. This can be deduced from the full theory of the growth form of periodic structures (first theorem in this paper, reformulated for graphs).
\item Moreover, graph vertices with distance n from the initial vertex lie exactly on the growth hexagon stretched n times.
\item Each tile in $HAT_1$ contains exactly one graph vertex. The tiles in $HAT_1$ are neighbours if and only if the corresponding vertices are neighbours in the graph.
\item Hence, the distance $d_{HAT_1}(T_1, T_2) = d_{Graph}(x_1, x_2)$, where $x_i$ is the graph vertex in $T_i$.
\item There exists $C$ such that any tile lies in the $C$-neighbourhood of the corresponding graph vertex.
\item The results follow from Figures 4 and 5.
\end{enumerate}

Geodetic chains are chains of minimal length between two tiles.
The described construction also provides some information about geodetic chains in $HAT_1$ because we can easily find many geodesic chains in the hexagonal tiling (or its dual graph).

\begin{figure}
\centering
\includegraphics[width=0.9\textwidth]{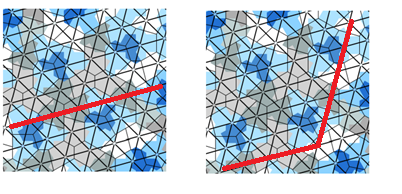}
\caption{Geodetic chains on the triangle tiling.}
\label{fig Hat proof 2}
\end{figure}

For example, sequential vertices lying on one straight line of the triangle graph form a geodetic chain. Moreover, vertices lying on a pair of rays with the angle $2\pi/3$ also form a geodetic chain, Figure~\ref{fig Hat proof 2}. The corresponding tiles are the geodetic chains in $HAT_1$.

Reitebuch studies lines of the triangle graph on the $HAT_1$ tiling. He considers two types of lines: lines nonintersecting blue tiles (black lines in his terminology) and lines intersecting blue tiles (blue lines).

\begin{figure}
\centering
\includegraphics[width=0.8\textwidth]{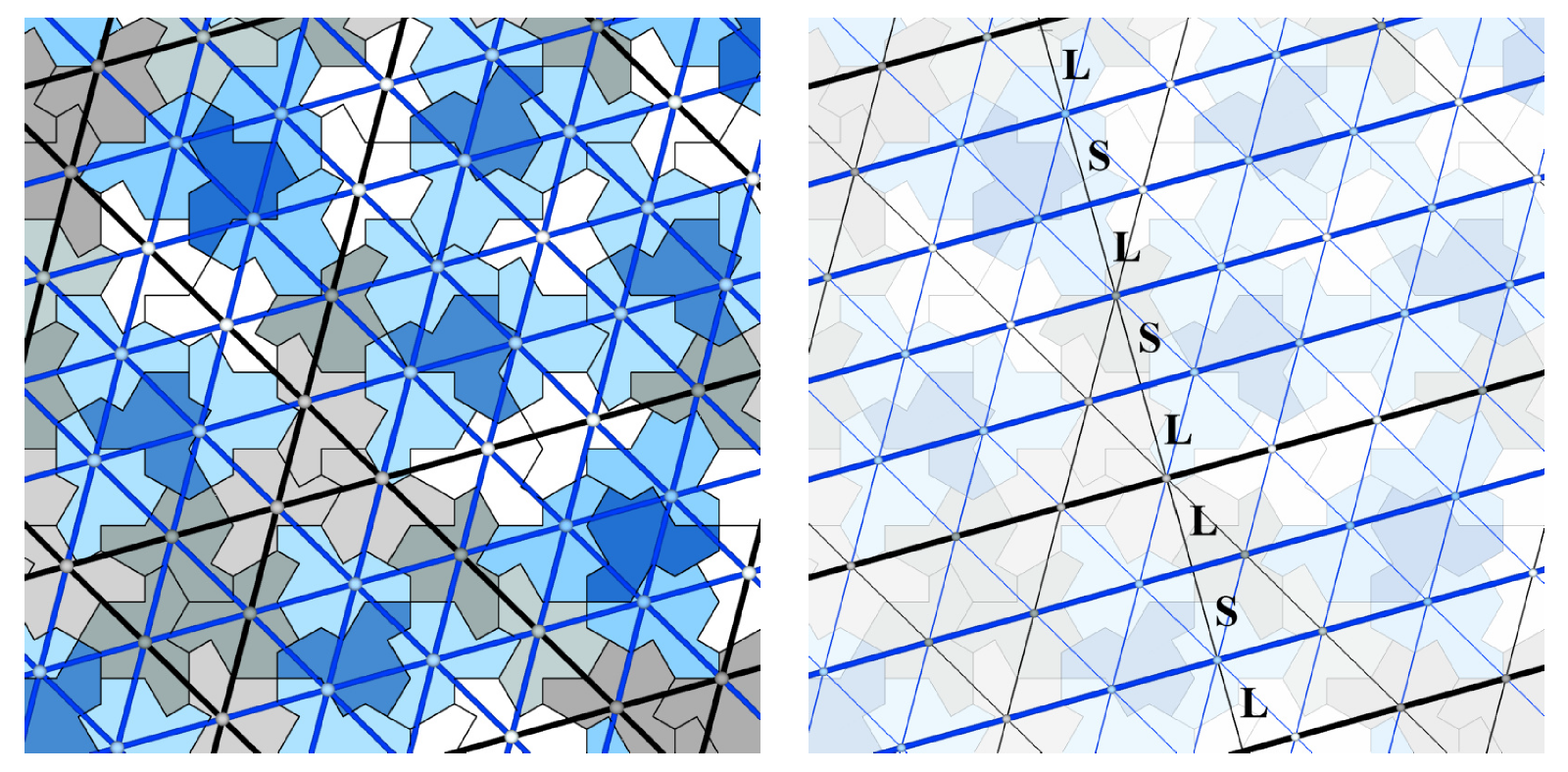}
\caption{Black and blue triangulation lines in the hat tiling -- labelled gaps between lines~\cite{Reitebuch2023}.}
\label{fig Hat proof 3}
\end{figure}

We call tiles of $HAT_1$ that intersect black lines important tiles. Note that important tiles are also tiles of $HAT$. Reitebuch proves that there are 2 or 4 blue lines between the black lines. He also proves that there are only 2 possible distances between parallel lines. Hence, some constant $K$ exists so that for any tile $T_1$ from $HAT_1$, there is an important tile $T_2$ such that $d_{HAT_1}(T_1, T_2) \leq K$. An analogue relation is also true for $HAT$ and Euclidean distances (but with another $K$).

Now, we need to explain why the growth forms of $HAT$ and $HAT_1$ coincide. For the proof, it is sufficient to show that for any important tiles $T_1$ and $T_2$, $d_{HAT}(T_1, T_2) = d_{HAT_1}(T_1, T_2)$. Sketch of the proof:

\begin{enumerate}
\item We can always assume that $T_1$ is important because the growth form does not depend on the initial tile.
\item If $T_2$ is important, we use that $d_{HAT}(T_1, T_2) = d_{HAT_1}(T_1, T_2)$, the growth form of $HAT_1$ is a hexagon and conclude that $T_2$ lies in the $C$-neighbourhood of the $n$ * growth hexagon. Here, $n=d_{HAT}(T_1, T_2)$.
\item Suppose that $T_2$ is a nonimportant tile. Denote by $T_{2\prime}$ some important tile such that
$d_{HAT}(T_2, T_{2\prime}) \leq K$. According to the above argument, such a tile always exists.
\item If $d_{HAT}(T_1, T_2)= n$ from the triangle inequality, we have $n-K \leq d_{HAT}(T_1, T_{2\prime}) \leq n + K$.
\item Hence, $T_{2\prime}$ lies in the $C$-neighbourhood of the $n^\prime$ * growth hexagon where $n-K\leq n^\prime \leq n + K$.
\item Hence, $T_{2\prime}$ lies in the $C^\prime$ neighbourhood of $n$ * growth hexagon with some new $C^\prime$.
\item As $d_{HAT}(T_2, T_{2\prime}) \leq K$, $T_2$ lies in the $C^{\prime\prime}$ neighbourhood of $n$ * growth hexagon with some new $C^{\prime\prime}$.
\end{enumerate}

For any two important tiles $T_1$ and $T_2$, there is a geodesic chain in $HAT_1$ corresponding to vertices lying on a pair of black rays with the angle $2\pi/3$.
Hence, all the tiles of this geodesic are non-blue.
Therefore, this chain is also a correct chain in $HAT$.
Therefore, we have $d_{HAT_1}(T_1, T_2) \geq d_{HAT}(T_1, T_2)$.

On the other hand, consider some chains $CH$ from $T_1$ to $T_2$ in $HAT$. We want to show that there exists a chain $CH_1$ from $T_1$ to $T_2$ in $HAT_1$ such that $d(CH_1) \leq d(CH)$.
This implies that $d_{HAT_1}(T_1,T_2)\leq d_{HAT}(T_1,T_2)$. Therefore, the distances $d_{HAT_1}(T_1, T_2)$ and $d_{HAT}(T_1, T_2$) are equal, and the growth forms coincide.

Consider patches (clusters) in $HAT$ consisting of tiles neighbouring a dark tile.
It is sufficient to consider the case in which $T_1$ and $T_2$ are tiles from one such patch.

Sketch of the explanation; why it is true:
\begin{enumerate}
\item Consider tiles $T_1$ and $T_2$ from $HAT$. As we showed before, it is sufficient to consider the case in which $T_1$ and $T_2$ are non-blue.
\item Consider a geodesic chain $CH$ from $T_1$ to $T_2$ in $HAT$. We want to construct a chain $CH_1$ in $HAT_1$ having the same length.
\item If $CH$ does not contain blue tiles, take $CH_1=CH$.
\item If $CH$ contains light blue tiles but not dark blue tiles, note that any light blue tile from $HAT$ intersects with exactly one light blue tile from $HAT_1$. The light blue tiles are changed from $HAT$ to the corresponding light blue tiles from $HAT_1$. Note that this change transforms chain to chain, that is, in the transformed sequence, consecutive tiles are neighbours.
\item If $CH$ contains dark blue tiles, the following transformation is used: Take a fragment of the geodesic chain containing all the dark blue tiles in one patch. Add to it tile $T_3$ ahead of the first in the fragment and tile $T_4$ next (in the geodesic) to the last tile in the fragment. The chain is changed from $T_3$ to $T_4$ to geodesic in the $HAT_1$ chain between such tiles. Our assumption is that $d_{HAT_1}(T_3, T_4) \leq d_{HAT}(T_3, T_4)$. Hence, after all such transformations, we obtain the chain $CH_1$ in $HAT_1$ having the same or less length as $CH$.
\end{enumerate}

\begin{figure}[H]
\centering
\includegraphics[width=0.9\textwidth]{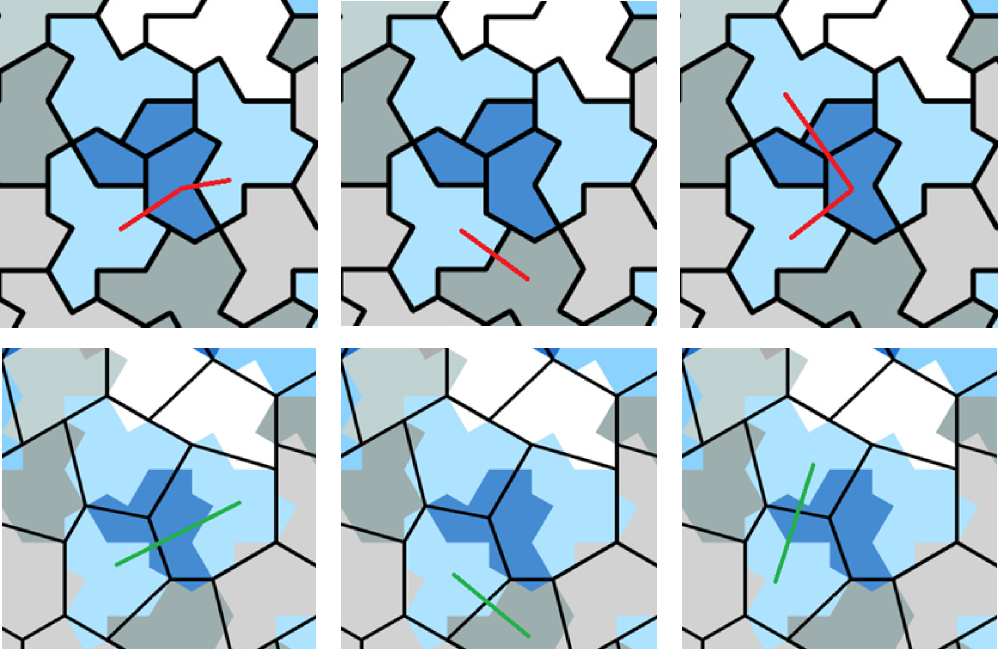}
\caption{top row: a cluster in $HAT$; bottom row: corresponding cluster in $HAT_1$.}
\label{fig Hat proof 4}
\end{figure}

In $HAT$, there is only a finite number of clusters (consisting of tiles neighbouring dark tiles) up to rigid motion.
Therefore, we need to verify the relation $d_{HAT}(T_1, T_2) = d_{HAT_1}(T_1, T_2)$ for a finite number of pairs $T_1, T_2$.
This verification can be performed directly by hand in the following way:

Draw some concrete dark blue tiles  in $HAT$, and all tiles neighbouring it.
Draw similar picture for $HAT_1$. For all pairs of tiles neighbouring clusters, find the geodesic chains in $HAT$ and $HAT_1$.
Verify that the lengths of the geodesics in $HAT_1$ are less than or equal to the corresponding lengths in $HAT$.
Figure~\ref{fig Hat proof 4} shows such verification for one concrete cluster.

\end{appendices}

\bibliography{growthforms.bib}

\end{document}